\documentstyle[12pt]{article}
\newlength{\abstractwidth}
\renewcommand{\thefootnote}{\fnsymbol{footnote}}
\renewcommand{\thanks}[1]{\footnote{#1}} 
\newcommand{\starttext}{ \setcounter{footnote}{0}
\renewcommand{\thefootnote}{\arabic{footnote}}}

\newcommand{\be}{\begin{equation}}
\newcommand{\bea}{\begin{eqnarray}}
\newcommand{\eea}{\end{eqnarray}} \newcommand{\ee}{\end{equation}}
 
 \def\ba{\begin{eqnarray}}
\def\ea{\end{eqnarray}}

\def\cM{{\cal M}}

\def\ra{\rightarrow}

\def\o{\omega}

\def\log{\,{\rm log}\,}

\def\o{\omega}

\def\d{\delta}
\def\e{\varepsilon}

\def\l{\lambda}
\def\m{\mu}

\def\o{\omega}

\def\na{\nabla}

\def\ve{\varepsilon}
\def\ge{\geq}
\def\le{\leq}
\def\ddt{{\partial\over\partial t}}
\def\ov{\overline}
\def\ti{\tilde}

\def\i{\infty}
\def\I{\int}
\def\p{\prod}
\def\s{\sum}

\def\ddb{{\partial\bar\partial}}

\def\ra{\rightarrow}

\def\cM{{\cal M}}

\def\na{{\nabla}}

 \def\v{\vskip .1in}

\def\[{{\bf [}}
\def\]{{\bf ]}}

\def\pl{\partial}

\begin{document}
\starttext \baselineskip=18pt \setcounter{footnote}{0}
\newtheorem{theorem}{Theorem}
\newtheorem{lemma}{Lemma}
\newtheorem{definition}{Definition}
\begin{center}
{\Large \bf THE K\"AHLER-RICCI FLOW AND THE }

{\Large \bf $\bar\partial$ OPERATOR ON VECTOR FIELDS\footnote{Research supported in part by
National Science Foundation grants  DMS-02-45371, DMS-06-04805, DMS-05-14003,
and DMS-05-04285.  Part of this work was carried out while the second-named author was supported by MSRI as a postdoctoral fellow. }}
\\
\bigskip
\bigskip

{\large D.H. Phong$^*$, Jian Song$^{**}$, Jacob Sturm$^\dagger$ and
Ben Weinkove$^\ddagger$} \\

\bigskip

\begin{abstract}
{\small The limiting behavior of the normalized K\"ahler-Ricci flow for manifolds with
positive
first Chern class is examined under
certain stability conditions. First, it is shown that if the Mabuchi K-energy is
bounded from below, then the scalar curvature converges uniformly to a constant.
Second, it is shown that if the Mabuchi K-energy is bounded from below and
if the lowest positive eigenvalue of the $\bar\partial^\dagger \bar\partial$ operator on smooth vector
fields is bounded away from 0 along the flow, then the metrics converge exponentially fast
in $C^\infty$ to a K\"ahler-Einstein metric.}

\end{abstract}

\end{center}

\baselineskip=15pt
\setcounter{equation}{0}
\setcounter{footnote}{0}

\section{Introduction}
\setcounter{equation}{0}

The K\"ahler-Ricci flow is arguably the most natural non-linear heat flow on a K\"ahler
manifold, and its singularities and asymptotic behavior can be expected to provide a
particularly deep probe of the geometry of the underlying manifold.
For manifolds $X$ with positive first Chern class, the K\"ahler-Ricci flow
exists for all times \cite{C}, and the issue is its asymptotic behavior. The convergence
of the flow would produce a K\"ahler-Einstein metric, the existence of which
had been conjectured by Yau \cite{Y2} to be equivalent to
the stability of $X$ in the sense
of geometric invariant theory. Thus the convergence and, more generally,
the asymptotic behavior of the flow should be related to stability conditions.

\v
There have been however only relatively few results in this
direction. In fact, the convergence of the flow for $c_1(X)>0$ has
been established only for $X={\bf CP}^1$ \cite{H, Ch, CLT}, for $X$
admitting a metric with positive bisectional curvature (and hence
must be ${\bf CP}^n$) \cite{CT}, under the assumption that $X$
already admits a K\"ahler-Einstein metric \cite{P2} or a
K\"ahler-Ricci soliton \cite{TZ2},
and for $X$ toric with vanishing Futaki invariant \cite{Z}
(which is known to imply that $X$ admits a K\"ahler-Einstein
metric \cite{WZ}). In \cite{PS}, it was shown that
certain stability conditions do imply the convergence of the flow,
without an a priori assumption on the existence of a
K\"ahler-Einstein metric or a K\"ahler-Ricci soliton, but with an
additional assumption on curvature bounds.

\v
The purpose of this paper is to relate the asymptotic behavior of the K\"ahler-Ricci flow
to stability conditions, without either of the previous assumptions of curvature bounds or
existence of K\"ahler-Einstein metrics or K\"ahler-Ricci solitons. More specifically, we deal with
two types of stability conditions. The first is the familiar lower boundedness of the
Mabuchi K-energy \cite{B, T, D2, PS}. The second is the
lower boundedness of the
first positive eigenvalue of the $\overline{\partial}^{\dagger} \bar\partial$ operator on smooth $T^{1,0}$ vector fields.
This second condition
appears to be new, but it should be closely related to the stability condition (B) introduced in
\cite{PS},
namely that the closure of the orbit of the complex structure $J$ of $X$ does not contain
any complex structure $\tilde J$ with a strictly higher number of independent holomorphic vector
fields.

\v
Our results are of two types. To state them precisely, let $X$ be a compact K\"ahler manifold
of dimension $n$ with $c_1(X)>0$, and let the
K\"ahler-Ricci flow\footnote{In this paper, we consider only the normalized
K\"ahler-Ricci flow, and designate this flow simply
by ``K\"ahler-Ricci flow".}  be defined by
\be
\label{kr}
{\pl\over\pl t}g_{\bar kj}=-(R_{\bar kj}-g_{\bar
kj}), \qquad g_{\bar kj}{\vert_{t=0}}=(g_0)_{\bar kj},
\ee
where $(g_0)_{\bar kj}$ is a given initial metric, with K\"ahler form
$\o_0={\sqrt{-1}\over 2}(g_0)_{\bar kj}dz^j\wedge d\bar z^k\in \pi\,c_1(X)$.
The Mabuchi K-energy $\cM(\o_\phi)=\cM(\phi)$ is the functional defined on $\pi c_1(X)$
by its value at some fixed reference point and its variation
\be
\label{mabuchi}
\delta \cM(\phi)=-{1\over V}\int_X \delta\phi (R-n)\o_\phi^n,
\qquad\qquad V\equiv\int_X \o_\phi^n=\pi^n c_1(X)^n,
\ee
where $\o_\phi=\o_0+{\sqrt{-1}\over 2}\ddb\phi\in\pi c_1(X)$ has been identified with its potential
$\phi$ (modulo constants)
and $R=g^{j\bar k}R_{\bar kj}$ denotes the scalar curvature of $\o_\phi$.
The first type of result assumes only a lower bound of the Mabuchi K-energy.
Under such an assumption, using the continuity method,
Bando \cite{B} had shown the existence of K\"ahler metrics
in $\pi c_1(X)$ with $||R-n||_{C^0}$ arbitrarily small. In \cite{PS}, \S 6,
it was shown that, under the same assumption, $||R-n||_{L^2}\ra 0$
for the K\"ahler-Ricci flow. Here we show:

\begin{theorem}
\label{almostcsck}
Assume that the Mabuchi K-energy is bounded from below on the K\"ahler
class $\pi c_1(X)$. Let $g_{\bar kj}(t)$ be any solution of the
K\"ahler-Ricci flow (\ref{kr}), and let $R(t)$ be the scalar curvature
of $g_{\bar kj}(t)$. Then we have
\begin{enumerate}
\item[(i)] $\displaystyle{\|R(t)-n\|_{C^0}\ra 0\ as\ t\ra\infty;}$
\item[(ii)] $\displaystyle{\int_0^\infty \|R(t)-n\|_{C^0}^p \, dt <\infty,
\qquad \textrm{when} \ p>2.}$
\end{enumerate}
\end{theorem}

\bigskip

In view of Lemma \ref{integrability} below,
the gap between a lower bound for the Mabuchi K-energy and
the existence of a K\"ahler-Einstein metric is thus
at most the gap between $L^p[0,\infty)$, $p>2$, and $L^1[0,\infty)$.
The second type of result assumes both a lower bound of the Mabuchi K-energy
and a stability condition (S):

\v

\begin{theorem}
\label{csck} Fix $\omega_0\in \pi c_1(X)$.
Let $g_{\bar kj}(t)$ be the solution of the K\"ahler-Ricci flow
with initial value $(g_0)_{\bar kj}$, and $\o(t)$ the corresponding
K\"ahler forms.
Let $\lambda_\omega$ be the lowest strictly positive
eigenvalue of the Laplacian $\bar\partial^\dagger\bar\partial
=-g^{j\bar k}\nabla_j\nabla_{\bar k}$ acting on smooth $T^{1,0}$ vector
fields.

\medskip

(i) If the following two conditions are satisfied,
\bea
&&\emph{(A)}\qquad\qquad {\rm inf}_{\omega \in \pi c_1(X)}\cM(\omega)> -\infty
\nonumber\\
&&\emph{(S)}\qquad\qquad
{\rm inf}_{t\in [0,\infty)}\lambda_{\omega(t)}>0
\nonumber
\eea
then the metrics $g_{\bar kj}(t)$ converge exponentially fast in $C^\infty$ to
a K\"ahler-Einstein metric.

(ii) Conversely, if the metrics $g_{\bar kj}(t)$ converge in $C^\infty$ to
a K\"ahler-Einstein metric, then the conditions {\rm (A)} and {\rm (S)} are satisfied.

(iii) In particular, if the metrics $g_{\bar kj}(t)$ converge in $C^\infty$
to a K\"ahler-Einstein metric,
then they converge exponentially fast in $C^\infty$ to this metric.

\end{theorem}

\v
As an immediate consequence, if the Mabuchi K-energy is bounded below on $\pi c_1(X)$ and 
\bea
&&
{\rm inf}_{\o\in\pi c_1(X)}\lambda_{\omega}>0
\nonumber
\eea
then Theorem \ref{csck} implies that every solution $g_{\bar kj}(t)$ of the K\"ahler-Ricci flow
converges exponentially fast in $C^\infty$ to a K\"ahler-Einstein metric.

\bigskip
We conclude this introduction with some remarks on the proof.
The recent works of Perelman \cite{P2} have provided powerful tools for
the study of the K\"ahler-Ricci flow, including a non-collapsing
theorem and the uniform boundedness
of the Ricci potential and of the scalar curvature.
On the other hand, there are still no known uniform bounds for the Riemannian and the
Ricci curvatures. We bypass this difficulty by exploiting two features of
the flow: the first is its parabolicity, so that certain stronger norms for
the key geometric quantities can be controlled by weaker norms
at an {\it earlier} time (e.g. Lemma \ref{lemmasmoothing}); and the second is that such bounds
at earlier times can still produce the desired convergence statements
when combined with suitable {\it differential-difference inequalities}
(see e.g. the inequality (\ref{weaker}) below).

\section{Perelman's results}
\setcounter{equation}{0}

Perelman \cite{P1}, \cite{P2} proved the following estimates
 for a solution of (\ref{kr})(see
\cite{ST} for a detailed exposition).
The first is bounds for the Ricci potential
$u=u(t)$ defined by
\be
R_{\bar kj}-g_{\bar kj}=-\pl_j\pl_{\bar k}u,
\qquad
{1\over V}\int_X e^{-u}\o^n=1,
\ee
and the second is a non-collapsing theorem:

\begin{enumerate}
\item[(i)] There exists a constant $C_0$ depending only on
$g_{\bar kj}(0)$ such that
\be
\| u \|_{C^0}+ \| \nabla u\|_{C^0}
+\|R\|_{C^0} \le C_0.
\ee

\item[(ii)] Let $\rho>0$ be given.  Then there exists $c>0$ depending only on
$g_{\bar kj}(0)$ and $\rho$ such that for all points $x \in X$, all times
$t\ge0$ and all $r$ with $0 < r \le \rho$, we have
\be
\int_{B_r(x)} \omega^n > c \, r^{2n},
\ee
where $B_r(x)$ is the geodesic ball of radius $r$ centered at $x$
with respect to the metric $g=g(t)$.

\end{enumerate}

For the reader's convenience, we note that the exact statement of
(ii) can be derived from (i) as follows.  Make a change of
variable $t= -\log(1-2s)$ and define a Riemannian metric $h=h(s)$ by
$h(s) = (1-2s) g(t(s))$.  Then $h(s)$ is a solution of Hamilton's
Ricci flow  for $s \in [0, 1/2)$, and using the scalar curvature
bound of (i), one can apply Theorem 8.3.1 of
\cite{To}, or the arguments contained in \cite{ST}.

\section{A smoothing lemma}
\setcounter{equation}{0}

The important idea of a smoothing lemma, exploiting the parabolicity
of the K\"ahler-Ricci flow, is due to Bando \cite{B}.
For our purposes, we need the version below, the key feature of
which is the fact that it does not require a lower bound
on the Ricci curvature:

\begin{lemma}
\label{lemmasmoothing}
There exist positive constants $\delta$ and $K$ depending only on
$n$ with the following property. For any $\varepsilon$ with $0 < \varepsilon \leq
\delta$ and any $t_0 \geq 0$, if
$$ \|u(t_0)\|_{C^0}   \le \varepsilon,$$
then
$$\|\na u(t_0+2)\|_{C^0}+\|R(t_0+2) - n \|_{C^0} \le K\ve.$$
\end{lemma}

\bigskip
\noindent
{\it Proof.} By making a translation in time we can assume, without loss of generality, that $t_0=0$.
It is well-known that $u$ evolves by
\begin{equation}
\label{eqnevolveu}
\ddt{u} = \Delta u + u - b,
\end{equation}
where $b=b(t)$ is the average of $u$ with respect to the measure
$e^{-u}\omega^n$:
\begin{equation}
\label{eqnb}
b = \frac{1}{V} \int_X u e^{-u} \omega^n.
\end{equation}

It is convenient to define a new constant $c=c(t)$ for $t\ge 0$ by
$\dot{c} = b+c$, $c(0) =0$.  Then set $\hat{u}(t) = - u(t) - c(t)$.
We have $\| \hat{u}(0)\|_{C^0} \le \ve$ and $\hat{u}$
evolves by
\begin{equation} \label{eqnevolveuhat}
\ddt{\hat{u}} = \Delta \hat{u} + \hat{u}.
\end{equation}

Following \cite{B}, we calculate
\begin{eqnarray} \label{eqnevolveu2}
&&\ddt{} \hat{u}^2  =  \Delta \hat{u}^2 - 2 | \nabla \hat{u}|^2 + 2
\hat{u}^2 \\ \label{eqnevolvenablau} &&\ddt{} | \nabla \hat{u} |^2
=  \Delta | \nabla \hat{u} |^2 - | \nabla \ov{\nabla} \hat{u} |^2 -
| \nabla \nabla \hat{u} |^2 + | \nabla \hat{u}|^2 \\
\label{eqnevolvelapu} && \ddt{} \Delta \hat{u}  =  \Delta ( \Delta
\hat{u} ) + \Delta \hat{u} + | \nabla \ov{\nabla} \hat{u} |^2.
\end{eqnarray}
Then we see from (\ref{eqnevolveu2}) that $\| \hat{u}(t) \|_{C^0} \le
e^2 \ve$ for $t \in [0,2]$. From (\ref{eqnevolvenablau}) we have
\be
\ddt{} \left(e^{-2t} (\hat{u}^2 + t | \nabla \hat{u}|^2)\right) \le \Delta \left( e^{-2t} (\hat{u}^2 + t | \nabla \hat{u}|^2) \right),
\ee
giving $\| \nabla \hat{u} \|_{C^0}^2 (t) \le e^4 \ve^2$ for $t\in
[1,2]$.

We will now prove a lower bound for $\Delta \hat{u}$.  Set $$H =
e^{-(t-1)} ( |\nabla \hat{u}|^2 - \ve n^{-1} (t-1) \Delta \hat{u}
)$$ and compute using (\ref{eqnevolvenablau}) and
(\ref{eqnevolvelapu}),
\begin{eqnarray}
\ddt{} H & = & \Delta H - e^{-(t-1)} \left( \ve n^{-1} \Delta
\hat{u} + (1+\ve n^{-1} (t-1)) | \nabla \ov{\nabla} \hat{u} |^2 + |
\nabla \nabla \hat{u} |^2 \right).
\end{eqnarray}
For $t \in [1,2]$, using the inequality $( \Delta \hat{u})^2 \le n |
\nabla \ov{\nabla} \hat{u}|^2$ we obtain
\begin{eqnarray}
\ddt{} H & \le & \Delta H + e^{-(t-1)} n^{-1} (-\Delta \hat{u})( \ve
+ \Delta \hat{u}).
\end{eqnarray}
We claim that $H < 2e^4 \ve^2$ for $t \in [1,2]$.  Otherwise, at the
point $(x', t') \in X \times (1,2]$ when this inequality first fails
we have $- \Delta \hat{u} \ge e^4 \ve$.  But since $(\ddt{} -
\Delta) H  \ge 0$ at this point, we also have $\ve + \Delta \hat{u}
\ge 0$, which gives a contradiction.   Hence at $t=2$ we have $H < 2
e^4 \ve^2$ and
$$\Delta \hat u > - 2 n e^5 \ve, $$
on $X$.

By considering the quantity $$K = e^{-(t-1)} ( |\nabla \hat{u}|^2 +
\ve n^{-1} (t-1) \Delta \hat{u} ),$$ we can similarly prove that
$\Delta \hat{u} < 2 n e^5 \ve$ at $t=2$ (see \cite{B}).  Since
$\Delta \hat{u} = R -n$ this completes the proof of the lemma.  Q.E.D.

\pagebreak[2]
\bigskip
\noindent
{\bf Remarks}

\begin{enumerate}
\item[(i)]  In the statement of the lemma, $(t_0+2)$ could be replaced by $(t_0+\zeta)$ for any
positive constant $\zeta$, at the expense of allowing the constants to depend on $\zeta$.
\item[(ii)]  Bando gives a different argument for the lower bound of $\Delta \hat{u}$ making use
of the fact that in his application there is a lower bound on the Ricci curvature at the initial time.
\item[(iii)]  A similar smoothing argument is also used in the proof of the Moser-Trudinger
inequality for K\"ahler-Einstein manifolds \cite{T}, \cite{TZ1}, \cite{PSSW}.
\end{enumerate}

\section{Proof of Theorem 1, part (i)}
\setcounter{equation}{0}

We provide now the proof of Theorem \ref{almostcsck}, part (i). In view of Lemma \ref{lemmasmoothing},
it suffices to show that $\|u\|_{C^0}\ra 0$ as $t\ra\infty$.
Recall that
$b$ is the average of $u$ with respect to the measure $e^{-u}\o^n$.
It suffices to show that $b$ and $\|u-b\|_{C^0}$ tend to $0$ as $t\ra\infty$.
This is an immediate consequence of Lemmas \ref{b} and \ref{Y} below
(together with Perelman's uniform bound for $||\na u||_{C^0}$), so
it suffices to prove these lemmas\footnote{It has been shown by H. Li \cite{L} that $b(t_m)\ra 0$
for a sequence of times $t_m\ra\infty$.}.

\v
We shall need the following Poincar\'e-type inequality for manifolds with $c_1(X)>0$
(see \cite{F} or \cite{TZ2}).

\bigskip

\begin{lemma}
\label{poincare}
Let $u$ satisfy the equation $R_{\bar kj}-g_{\bar kj}=-\pl_j\pl_{\bar k}u$.
Then the following inequality
\be
{1\over V}\int_X f^2 e^{-u}\o^n\leq
{1\over V}\int_X|\na f|^2 e^{-u}\o^n
+
({1\over V}\int_X f e^{-u}\o^n)^2,
\ee
holds for all $f\in C^\infty(X)$.
\end{lemma}

\v
\noindent
{\it Proof of Lemma \ref{poincare}}: We include the easy proof for the reader's convenience.
The desired inequality is equivalent to the fact that the lowest strictly positive eigenvalue
$\mu$ of the following operator
\be
-g^{j\bar k}\na_j\na_{\bar k}f+g^{j\bar k}\na_{\bar k}f\na_ju =\mu f,
\ee
with eigenfunction $f$ satisfies $\mu \geq 1$. (Note that this operator
is self-adjoint with respect to the measure $V^{-1}e^{-u}\o^n$, and that its kernel
consists of constants.)
Applying $\na_{\bar l}$ and commuting $\na_{\bar l}$ through in the first term gives
\be
-g^{j\bar k}\na_j\na_{\bar k}\na_{\bar l}f
+
R_{\bar l}{}^{\bar p}\na_{\bar p}f+
g^{j\bar k}\na_{\bar l}\na_j u\na_{\bar k}f+g^{j\bar k}\na_{\bar l}\na_{\bar k}f \na_j u
=
\mu\na_{\bar l}f.
\ee
Integrate now against $g^{m\bar l}\na_m f\,e^{-u}\o^n$ and integrate by parts.
In view of the fact that $R_{\bar k j}+\pl_j\pl_{\bar k}u=g_{\bar kj}$, we obtain
\be
\int_X|\bar\na\bar\na f|^2 e^{-u}\o^n+\int_X|\bar\na f|^2 e^{-u}\o^n
=\mu\int_X|\bar\na f|^2e^{-u}\o^n,
\ee
from which the desired inequality $\mu\geq 1$ follows at once. Q.E.D.

\bigskip
Henceforth we shall denote by $\|\cdot\|_{L^2}$ the $L^2$ norm with
respect to the measure $\o^n$. This $L^2$ norm is uniformly equivalent to
the $L^2$ norm with respect to the measure $e^{-u}\o^n$, in view
of Perelman's theorem. The following lemma holds in all generality for
the K\"ahler-Ricci flow:

\bigskip

\begin{lemma}
The Ricci potential $u=u(t)$ and its average $b=b(t)$ satisfy
the following inequalities, where the constant $C$ depends only on $g_{\bar kj}(0)$:
\label{b}
\begin{enumerate}
\item[(i)] $\displaystyle{0 \le -b \le \|u-b\|_{C^0}}$;
\item[(ii)] $\displaystyle{\|u-b\|_{C^0}^{n+1}\ \leq\ C\,\|\nabla u\|_{L^2}\,\|\nabla u\|_{C^0}^n}$.
\end{enumerate}
\end{lemma}

\bigskip
\noindent
{\it Proof of Lemma \ref{b}}: First, we observe that, as a consequence of
Jensen's inequality and the convexity of the exponential function,
\be
b={1\over V}\int_X u\,e^{-u}\o^n
\leq
\log({1\over V}\int_X e^u e^{-u}\o^n)=0.
\ee
On the other hand, $e^{-u}$ has average $1$ with respect to the measure $\o^n$,
and thus ${\rm sup}_X u\geq 0$. Thus $- b \leq {\rm sup}_X(u-b)$,
and (i) is proved.

\v

Next, let $A=\| u-b \|_{C^0}=|u-b|(x_0)$. Then $|u-b|\geq {A\over 2}$ on the ball $B_r(x_0)$ of radius
$r={A\over 2\|\na u\|_{C^0}}$ centered at $x_0$. If $r<\rho$, where $\rho$
is some fixed uniform radius in Perelman's non-collapsing result, then
\be
\I_X(u-b)^2\o^n\geq \I_{B_r(x_0)} {A^2\over 4} \o^n \, \geq\ c{A^2\over 4}
\left({A\over 2\|\na u\|_{C^0}}\right)^{2n}
\ee
and thus
\be
\|u-b\|_{C^0}^{n+1}\ \leq \ C_1\|\na u\|_{C^0}^n\|u-b\|_{L^2}.
\ee
Applying Lemma \ref{poincare}, we have
\be
\label{4.1}
\|u-b\|_{C^0}^{n+1}
\leq \ C_1\|\na u\|_{C^0}^n\|u-b\|_{L^2}\ \leq \
C_2\|\na u\|^n_{C^0}\|\na u\|_{L^2}.
\ee
On the other hand, if $r>\rho$, then integrating over the ball $B_{\rho}(x_0)$
gives the bound $\|u-b\|_{C^0}\leq C\,\|\na u\|_{L^2}$, which is a stronger
estimate than the one we need. Q.E.D.

\bigskip

\begin{lemma}
\label{Y}
Assume the Mabuchi K-energy is bounded from below on $\pi c_1(X)$.
Set
\be
Y(t)=\int_X |\nabla u|^2\o^n=\|\na u\|_{L^2}^2.
\ee
Then for any choice of initial condition $\o_0\in\pi c_1(X)$,
the quantity $Y(t)\ra 0$ along the K\"ahler-Ricci flow as $t\to\infty$.
\end{lemma}

\v
\noindent
{\it Proof of Lemma \ref{Y}}: The proof of this Lemma can be found in \cite{PS},
\S 6. We provide the short proof, for the sake of completeness.
Let $\phi=\phi(t)$ be the potential of $g_{\bar kj}(t)$, so that
$g_{\bar kj}=(g_0)_{\bar kj}+\pl_j\pl_{\bar k}\phi$. Then clearly
$\dot\phi-u$
is a constant depending only on time along the K\"ahler-Ricci
flow, and it follows immediately from the definition of the Mabuchi K-energy
(\ref{mabuchi}) that its
derivative along the K\"ahler-Ricci flow is given by
\be
{d\over dt}\cM (\phi)=-{1\over V}\int_X |\na u|^2\o^n=-{1\over V}Y(t).
\ee
Thus, since $\cM$ is bounded from below, we have for all $T>0$
\be
\label{Yintegral}
{1\over V}\int_0^T Y(t)dt=\cM(\phi_0)-\cM(\phi)\leq C,
\ee
and hence $Y(t)$ is integrable over $[0,\infty)$. Equivalently
$\sum_{m=0}^\infty \int_m^{m+1}Y(t)dt<\infty$, and hence there
exists $t_m\in [m,m+1)$ with $Y(t_m)\to 0$. Next, $Y$ satisfies the
following differential identity (see \cite{PS}, eq. (2.10))
\be
\dot Y
=
(n+1)Y-\int_X|\na u|^2 R \, \o^n
-
\int_X|\bar\na\na u|^2\o^n-
\int_X|\na\na u|^2\o^n
\leq C\,Y
\ee
for some constant $C$, since $|R|$ is uniformly bounded by Perelman's estimate.
This implies $Y(t)\leq Y(s) \,e^{C(t-s)}$ for $t\geq s$. In particular
$Y(t)\leq Y(t_m) e^{2C}$ for all $t\in[m+1,m+2)$, and hence $Y(t)\ra 0$
as $t\ra\infty$. The proof of Lemma \ref{Y} and hence  of Theorem \ref{almostcsck}, part (i) 
is complete.  Q.E.D.

\section{Proof of Theorem 1, part (ii) and of Theorem 2}
\setcounter{equation}{0}

We begin by establishing Theorem \ref{csck}, part (i), that is,
the statement that the stability conditions (A) and (S) together
imply the exponential convergence in $C^\infty$ of $g_{\bar kj}(t)$
to a K\"ahler-Einstein metric.
This is an immediate consequence of
Lemmas \ref{exponential} and \ref{integrability} below, so it suffices to establish
those two lemmas.

\v

\begin{lemma}
\label{exponential}
Assume the Mabuchi K-energy is bounded
from below on $\pi c_1(X)$ and  we have $\lambda_t \ge \lambda>0$ along the K\"ahler-Ricci flow with initial value $(g_0)_{\bar kj}$.
Then,
the quantity $Y(t)=||\nabla u||_{L^2}^2$ tends to $0$ exponentially, that is,
there exist constants $\mu>0$ and $C>0$ independent of $t$ so that
\be
Y(t)\ \leq \ C\,e^{-\mu t},
\qquad t\in [0,\infty).
\ee
Moreover,
\be
\|u\|_{C^0}+\|\na u\|_{C^0}+\|R-n\|_{C^0}
\leq\ C\,e^{-{1\over 2(n+1)}\mu t},
\qquad t\in[0,\infty).
\ee
\end{lemma}

\bigskip
\noindent
{\it Proof of Lemma \ref{exponential}}: We recall the following inequality
from \cite{PS}, which holds for the K\"ahler-Ricci flow
without any additional assumption:
\be \dot Y\ \leq \ -2\l_tY-2\l_t{\rm Fut}(\pi_t(\na^ju))-\I_X|\na u|^2(R-n)\o^n-\I_X\na^ju\na^{\bar k}
u(R_{\bar k j}-g_{\bar k j})\o^n.
\ee
Here ${\rm Fut}(\pi_t(\na^j u))$ is the Futaki invariant, applied to the
orthogonal projection $\pi_t(\na^j u)$
of the vector field $\na^j u$ on the space of holomorphic vector
fields.

Now assume that the Mabuchi K-energy
is bounded below. Then the Futaki invariant
${\rm Fut}$ is identically $0$, and in view of Theorem \ref{almostcsck}, part (i), 
the preceding inequality reduces to, for $t$ sufficiently large,
\be
\dot Y\ \leq \ -\l \,Y-\I_X\na^ju\na^{\bar k}
u(R_{\bar k j}-g_{\bar k j})\o^n.
\ee
To get exponential convergence, we would like to show
$|\I_X\na^ju\na^{\bar k}
u(R_{\bar k j}-g_{\bar k j})\o^n|\leq {\l\over 2}Y$.
This would of course follow if we knew that
$R_{\bar kj}-g_{\bar kj}$ were small. In the absence of
such information, we
shall prove something
a bit weaker, but which turns out to be sufficient for our purposes.
We claim that there
exists $K_0>0$ such that
\be
\label{weaker}
 \dot Y(t)\ \leq \ -\l Y(t)\ + \ {\l\over 2}Y^{1\over 2}(t)\cdot
 \prod_{j=1}^N[Y(t-a_j) ]^{\d_j\over 2}\ \ \hbox{for all
$t\geq K_0$},
\ee
where $N$ is an integer, the $a_j$
are non-negative integers and the $\d_j$ are non-negative real numbers with
the property that $\s_{j=1}^N\d_j=~1$.
\v
To see this, first note that
$|\na^ju\na^{\bar k}
u(R_{\bar k j}-g_{\bar k j})|\leq |\na u|^2|R_{\bar kj}-g_{\bar kj}|$, and thus
\be
\left|\I_X\na^ju\na^{\bar k}
u(R_{\bar k j}-g_{\bar k j})\o^n\right|\ \leq \
||\nabla u||_{C^0}
(\int_X|\na u|^2\o^n)^{1/2}
(\int_X|R_{\bar kj}-g_{\bar kj}|^2\o^n)^{1/2}.
\ee
However, an integration by parts shows readily that
\be
\int_X|R_{\bar kj}-g_{\bar kj}|^2\o^n
=
\int_X|\pl_j\pl_{\bar k} u|^2\o^n=
\int_X|\Delta u|^2\o^n
=
\int_X|R-n|^2\o^n,
\ee
and hence
\be
\label{iteration0}
\left|\I_X\na^ju\na^{\bar k}
u(R_{\bar k j}-g_{\bar k j})\o^n\right|\ \leq \
Y^{\frac{1}{2}}(t)\,\|\nabla u\|_{C^0}
\|R-n\|_{L^2}.
\ee
We use Lemmas \ref{lemmasmoothing}, \ref{b} and \ref{Y}
 to estimate $||R-n||_{L^2}$
by $||\na u||_{L^2}$ and $||\na u||_{C^0}$
at an earlier time $t-2$. More precisely, we have, for $t$ sufficiently large
\bea
\|R-n\|_{L^2}(t)
&\leq& \|R-n\|_{C^0}(t)
\leq K\|u\|_{C^0}(t-2)
\leq
2K\|u-b\|_{C^0}(t-2)\nonumber
\\
&\leq&
C\,\|\na u\|_{C^0}^{n\over n+1}(t-2)
\|\na u\|_{L^2}^{1\over n+1}(t-2),
\eea
so that the coefficient of $Y^{\frac{1}{2}}(t)$ on the right hand side of
(\ref{iteration0}) can be estimated by
\be
\|\nabla u\|_{C^0}
\|R-n\|_{L^2}
\leq
C\,\|\na u\|_{C^0}(t)
\|\na u\|_{C^0}^{n\over n+1}(t-2)
\|\na u\|_{L^2}^{1\over n+1}(t-2).
\ee
We  wish to iterate this estimate by
applying the following bound:
\be
\label{it} \|\na u\|_{C^0}(t)\ \leq \ C_1 \|u-b\|_{C^0}(t-2)\
\leq \
C_2\|\na u\|_{C^0}^{n\over n+1}(t-2)\|\na u\|^{1\over
n+1}_{L^2}(t-2).
\ee

Let $A(t)=\|\na u\|_{L^2}(t)$
and $B(t)=\|\na u\|_{C^0}(t)$.
Suppose $G(t)$ is a function
of the form

$$ G(t)=\prod_j A(t-a_j)^{\sigma_j}
\prod_kB(t-b_k)^{\e_k}
$$
where $a_j,b_k$ are non-negative
integers and $\sigma_j,\e_k$ are non-negative
real numbers. We let $\sigma=\s \sigma_j$ and $\e=\s \e_k$ and assume that $\sigma+\e=2$.
Then (\ref{it})
implies
$ G(t)\ \leq \ C\tilde G(t)
$
where $ \tilde G(t)=\prod_j A(t-\tilde a_j)^{\tilde \sigma_j}
\prod_kB(t- \tilde b_k)^{\tilde \e_k}
$
where $\ti \sigma =\s \ti\sigma_j$,
$\ti \e=\s_k\ti\e_k$ still have
the property $\ti\sigma+\ti\e=2$ but
$\ti\e= {n\over n+1}\e$.
Thus if we iterate, we see
that $G(t)\leq C\ti G(t)$ where
$\ti\e<1$ and $\ti\sigma>1$.
Setting $\d_j=\ti\sigma_j/\ti\sigma$, we get

\be \|\na u\|_{C^0}(t)\|\na u\|_{C^0}^{n\over n+1}(t-2)
\|\na u\|^{1\over
n+1}_{L^2}(t-2)\ \leq \
H(t)\cdot\prod_jA(t-\tilde a_j)^{\d_j}
\ee
where $H(t)=C\p_kB(t- \tilde b_k)^{\tilde{\e}_k}\p_jA(t- \tilde a_j)^{\ti\sigma_j-\d_j}$. Lemma \ref{Y} implies
$H(t)\ra 0$
as $t\to\i$, and thus we obtain (\ref{weaker}).
\v

Let $F(t)=Re^{-\m t}$ where $R>0$ and
$\m\in (0,1)$ are
positive numbers to be chosen later.
We want to show $Y\leq F$. Assume not.  Since $Y$ is bounded we may choose $R$ sufficiently large
so that for some time $t_0> K_0$ we have $Y(t) < F(t)$ for $0 \le t < t_0$ and $Y(t_0)=F(t_0)$.
\v
We claim
\be \dot Y(t_0)\leq -3\m Y(t_0). \label{eqnYt0} \ee
If not, then $\dot Y(t_0)\geq -3\m Y(t_0) = -3\m F(t_0)$ so that
$$ -3\m F(t_0)\
\leq \ -\l F(t_0)+{\l\over 2}F(t_0)^{1\over 2}\prod_{j=1}^N[Y(t_0-a_j) ]^{\d_j\over 2} \
\leq \ -\l F(t_0)+{\l\over 2}\prod_{j=0}^N[F(t_0-a_j) ]^{\d_j\over 2} \
$$
where we set $a_0=0$ and $\d_0=1$ so
that $\s_{j=0}^N\d_j=2$. Now we have

$$(\l-3\m)Re^{-\m t_0}\ \leq \ {\l\over 2}Re^{-\m t_0}e^{\m\s a_j\d_j/2}.\
$$
This implies
$$\  e^{-\m\s a_j\d_j/2}\ \leq \ {\l\over 2(\l-3\m)},
$$
and choosing $\m$ sufficiently close to zero gives a contradiction. This proves (\ref{eqnYt0}).
\v

On the other hand, from  the definition of $t_0$, we have
$$\left. \frac{d}{dt} \right|_{t=t_0} (Y-F) \ge 0,$$ 
and hence $\dot{Y} (t_0) \ge - \mu Y(t_0)$, which contradicts (\ref{eqnYt0}).
This proves $Y\leq F$ and thus
$Y$ decays exponentially. Now
Lemma \ref{b} implies that $||u||_{C^0}$
decays exponentially which, together
with Lemma \ref{lemmasmoothing}, shows that $||R-n||_{C^0}$
decays exponentially. The proof of Lemma \ref{exponential}
is complete. Q.E.D.

\begin{lemma}
\label{integrability}
Assume that the scalar curvature $R(t)$ along the K\"ahler-Ricci flow satisfies
\be
\int_0^\infty \|R(t)-n\|_{C^0}\,dt\ <\ \infty.
\ee
Then the metrics $g_{\bar k j}(t)$ converge exponentially fast in $C^{\infty}$ to a
K\"ahler-Einstein metric.
\end{lemma}

\bigskip
\noindent
{\it Proof of Lemma \ref{integrability}}:
The basic observation is that the integrability of $\|R-n\|_{C^0}$
over $t\in [0,\infty)$ implies a uniform bound for $\|\phi\|_{C^0}$,
where $\phi$ is the K\"ahler potential. More precisely, let the
potential $\phi(t)$ of $g_{\bar kj}(t)$
be normalized by
\be
{\pl\over \pl t}\phi=\log {\o^n\over\o_0^n}+\phi+u(0),
\qquad
\phi\vert_{t=0}=c_0,
\ee
with the constant $c_0$ chosen as in \cite{CT}, \cite{L}, and \cite{PSS}, eq. (2.10).
Then $g_{\bar kj}=(g_0)_{\bar kj}+\pl_j\pl_{\bar k}\phi$
satisfies the K\"ahler-Ricci flow, and Perelman's estimate
for $u$ implies $\|\dot\phi\|_{C^0}\leq C$.

\v

Now we have
\be
{d\over dt}\,(\log\, {\o^n\over\o_0^n})
=g^{j\bar k}\dot g_{\bar kj}=-(R-n).
\ee
Thus for any $t\in(0,\infty)$,
\be
\left| \log\,{\o^n\over\o_0^n} \right|
= \left| \int_0^t(R-n)\,dt \right| \leq \int_0^\infty \|R-n\|_{C^0} dt<\infty.
\ee
On the other hand, the K\"ahler-Ricci flow can be rewritten as
\be
\phi=-\log\,{\o^n\over\o_0^n}+\dot\phi -u(0)
\ee
and thus the uniform bound for $\|\phi\|_{C^0}$ follows from the
uniform bound for $|\log\,(\o^n/\o_0^n)|$ and Perelman's uniform
estimate for $\|\dot\phi\|_{C^0}$.

\v
The uniform boundedness of $\|\phi\|_{C^0}$ implies the uniform
boundedness of $\|\phi\|_{C^k}$ for each $k\in{\bf N}$
(see e.g. \cite{Y1, C, PSS, Pa}).   The metrics $g_{\bar kj}(t)$ are all uniformly
equivalent and bounded in $C^{\infty}$.
Thus there must
exist a subsequence of times $t_m\ra +\infty$ with $\phi(t_m)$ converging
in $C^\infty$ and the limit $\phi(\infty)$ is a potential for a smooth
K\"ahler-Einstein metric.   We claim that $\lambda_t$ is uniformly bounded below away from zero
along the flow.  If not, there would be a sequence of metrics $g_{\bar k j}(t_l)$ along the flow
with $\lambda_{t_l} \rightarrow 0$.  By the estimates above, after taking a subsequence,
the $g_{\bar k j}(t_{l})$ would converge in $C^{\infty}$ to a  K\"ahler metric $g'_{\bar k j}$
and $0= \lim_{l \rightarrow \infty} \lambda_{t_l} = \lambda (g'_{\bar k j}) >0$ (here one can apply the argument of \cite{PS}, \S 4 in the special case when the complex structure is fixed) giving a contradiction.  Moreover, the Mabuchi K-energy is bounded
below \cite{BM} and so we can apply Lemma \ref{exponential} to see that $Y(t)=\|\na u\|_{L^2}^2$
decays exponentially to $0$.   The arguments of \cite{PS}, \S 3 now show that $\|\nabla u\|_{(s)}$
decay exponentially to $0$
for any Sobolev norm $\|\cdot\|_{(s)}$,
and hence for any norm $\|\cdot\|_{C^k}$,
since  the metrics $g_{\bar kj}(t)$ are all
equivalent, and their Riemannian curvatures all uniformly bounded. But
then $\|\dot g_{\bar kj}\|_{C^k}=
\|R_{\bar kj}-g_{\bar kj}\|_{C^k}$ decays exponentially
to $0$ for any $k$, and hence the metrics $g_{\bar kj}$
converge exponentially fast to a K\"ahler-Einstein metric.  This completes the proof of
Lemma \ref{integrability} and hence of part (i) of Theorem \ref{csck}. Q.E.D.

\v

\bigskip
\noindent
{\it Proof of Theorem \ref{csck}, parts (ii) and (iii):}
\ Part (ii) follows immediately from the argument above.
Part (iii) follows immediately from Part (i) and Part (ii). Q.E.D.

\bigskip
\noindent
{\it Proof of Theorem \ref{almostcsck}, part (ii):}
\
First observe that
equation (\ref{it}) and Lemma \ref{lemmasmoothing}  imply
$$
\|R-n\|_{C^0}(t)\ \leq \ C_1
\|u\|_{C^0}(t-2)\ \leq 2C_1
\|u-b\|_{C^0}(t-2)
\ \leq \ C_2\|\na u\|_{C^0}^{n\over n+1}(t-2)\|\na u\|_{L^2}^{1\over
n+1}(t-2).
$$
Equation (\ref{it}), together with the iteration
argument used in the proof of Theorem \ref{csck},
show that if $t$ is sufficiently
large then
\be
\|R-n\|_{C^0}(t)\ \leq \ C \p_jA(t-a_j)^{\d_j}\p_kB(t-b_k)^{\e_k},
\ee
where $a_j,b_k$, $\d_j,\e_k$ are non-negative
real numbers, $\d+\e=\s_j\d_j+\s_k\e_k=1$ and $\d={2\over p}$. In
particular, if $N=\max_ja_j$ we have
\be
\|R-n\|_{C^0}(t)\ \leq \ C\p_j Y(t-a_j)^{\d_j\over 2}\ \
\hbox{for $t\geq N$}.
\ee
Since $\I_0^\i Y dt<\i$ if the Mabuchi K-energy is bounded below
we have

$$ \I_N^\i\|R-n\|_{C^0}^p\ dt\ \leq\ \
C_1\I_N^\i \p_jY(t-a_j)^{p\d_j\over 2}dt\ \leq \
C_1\p_j\left(\I_N^\i Y(t-a_j) dt \right)^{p\d_j\over 2}\ < \ \i.
$$
This establishes the desired inequality in Theorem 1, part (ii). Q.E.D.

\section{Further results and remarks}
\setcounter{equation}{0}

We conclude with some further results not required for Theorems
\ref{almostcsck} and \ref{csck} but which may be of interest. We also discuss
some possibilities for further developments.

\bigskip
(1) The average $b$ of the Ricci potential $u$ with respect to
the volume form $e^{-u}\o^n$ is monotone under the K\"ahler-Ricci flow\footnote{It has recently been brought to our attention
by V. Tosatti that this had also been observed earlier by Pali \cite{Pa}.},
\be
{d\over dt} b\geq 0.
\ee
To see this, we calculate using (\ref{eqnevolveu})
\bea
{d\over dt}b&=&{1\over V}{d\over dt}\int_X ue^{-u}\o^n\nonumber\\
&=&{1\over V}\int_X (\Delta u+u-b)e^{-u}\o^n-{1\over V}\int_X ue^{-u}(\Delta u+u-b)\o^n
+
{1\over V}\int_X ue^{-u}\Delta u\,\o^n
\nonumber\\
&=&
{1\over V}
\int_X |\na u|^2e^{-u}\o^n
-
{1\over V}\int_X u^2e^{-u}\o^n+b^2\geq 0,
\eea
where for the third line we have used the equality
\be
\int_X (\Delta u)e^{-u}\o^n=\int_X|\na u|^2 e^{-u}\o^n,
\ee
and for the last line we have used Lemma \ref{poincare}.


\bigskip
(2) It may be worth pointing out that Bando's result \cite{B},
or the $L^2$ result of \cite{PS}, \S 6, or Theorem \ref{almostcsck} proved here,
combined with Donaldson's recent result \cite{D2},
\be
\label{donaldson}
{\rm inf}_{\o\in \pi c_1(X)}\|R-n\|_{L^2}^2
\geq {\rm sup}_{\cal X} (-{1\over N^2_2({\cal X})}{\cal F}({\cal X}) )
\ee
each shows that the lower boundedness of the Mabuchi K-energy
in $\pi c_1(X)$ implies the K-semistability of
$X$. Here ${\cal X}$ denotes a test-configuration,
and ${\cal F}({\cal X})$ is its Futaki invariant (see \cite{D2} for the definition of  $N_2 ({\cal X})>0$). Indeed,
by the previously mentioned results, the lower boundedness of the K-energy
implies that the left hand side is $0$, and hence ${\cal F}({\cal X})\geq 0$
for any test configuration ${\cal X}$, which is the definition of K-semistability
\cite{T, D1}.

\bigskip
(3) It is natural to ask whether a converse to Bando's result (or to Theorem \ref{almostcsck})
is true, that is, whether the existence of metrics with $\|R-n\|_{C^0}\ra 0$
implies the lower boundedness of the Mabuchi K-energy. This would
complement well the result of
Tian \cite{T}, namely that the existence of a K\"ahler-Einstein
metric on a manifold $X$ with $c_1(X)>0$ and no holomorphic vector fields
is equivalent to the properness of the Mabuchi K-energy (in the case of
existence of holomorphic vector fields, there are some technical assumptions
on the automorphism group of $X$, see \cite{T}, and also \cite{PSSW}).

\bigskip
(4) The condition (S) in Theorem \ref{csck} can be interpreted as a stability
condition in the following sense, along the lines of \cite{PS}. Assume that there exist
diffeomorphisms $F_t:X\ra X$
so that $(F_t)_*(g(t))$ converges in $C^\infty$ to a metric
$\tilde g(\infty)$. Then if $J$ is the
complex structure of $X$, the pull-backs $(F_t)_*(J)$ converge also to a
complex structure $J(\infty)$ (see \cite{PS}, \S 4). Clearly, the eigenvalues
$\lambda_{\o(t)}$ are unchanged under $F_t$.
If they don't remain bounded away
from $0$ as $t\to\infty$, then the complex structure $J(\infty)$ would
have a strictly higher
number of independent vector fields than $J$. Thus the $C^\infty$ closure of the orbit
of $J$ under the diffeomorphism group contains a complex structure different from
$J$, and $J$ cannot be included in a Hausdorff moduli space of complex structures.

\bigskip

(5) One set of assumptions which would guarantee the existence of diffeomorphisms
$F_t$ is the uniform boundedness of the Riemannian curvature tensor along
the K\"ahler-Ricci flow. Under such an assumption, it was shown in \cite{PS} that
the condition (B) introduced there, and quoted earlier in our Introduction,
implies our condition (S).
Clearly, it would be very valuable to relate (B) and (S) in more
general situations.

\bigskip
(6) Under the sole condition of lower boundedness of the Mabuchi K-energy,
it was shown in \cite{PS}, \S 6, that there exists a sequence of times $t_m\ra\infty$
with
\be
\|\na R(t_m)\|_{L^2}\ra 0.
\ee
It would be interesting to determine whether this convergence
can take place with stronger norms.

\bigskip
(7)  We observe that the following alternative version of Theorem 2, part (iii) also holds by our results:   if a solution $g(t)$ of the K\"ahler-Ricci flow converges to a K\"ahler-Einstein metric in $C^{\infty}$ \emph{modulo automorphisms} (i.e. there exists a family of biholomorphisms $\Psi_t : X \rightarrow X$ such that $(\Psi_t)_*(g(t))$ converges in $C^{\infty}$ to a K\"ahler-Einstein metric) then the unmodified K\"ahler-Ricci flow $g(t)$  converges exponentially fast to a K\"ahler-Einstein metric.

\bigskip
\noindent
{\bf Acknowledgements} \ The authors thank the referee for some helpful suggestions.


\newpage

\begin{thebibliography}{99}

\bibitem[B]{B} Bando, S. {\em The K-energy map, almost Einstein K\"ahler metrics and an inequality of the
Miyaoka-Yau type}, T\^ohoku Math. J. {\bf 39} (1987), 231--235
\bibitem[BM]{BM} Bando, S. and Mabuchi, T. {\em Uniqueness of Einstein K\"ahler metrics modulo connected
group actions},  Algebraic geometry, Sendai, 1985,  11--40, Adv. Stud. Pure Math., 10, North-Holland,
Amsterdam, 1987
\bibitem[C]{C} Cao, H.-D. {\em Deformation of K\"ahler metrics to K\"ahler-Einstein metrics on compact
K\"ahler manifolds},  Invent. Math. {\bf 81} (1985),  no. 2, 359--372
\bibitem[CLT] {CLT} Chen, X.X., Lu, P. and Tian G. {\em A note on uniformization of Riemann surfaces by
Ricci flow}, Proc. Amer. Math. Soc. {\bf 134} (2006), no. 11, 3391--3393
\bibitem[CT] {CT} Chen, X.X. and Tian, G.
{\em Ricci flow on K\"ahler-Einstein manifolds}, Duke Math. J.  {\bf 131}  (2006),  no. 1, 17--73
\bibitem[Ch] {Ch} Chow, B. {\em The Ricci flow on the 2-sphere},
J. Differential Geom. {\bf 33} (1991) 325--334
\bibitem[D1]{D1} Donaldson, S.K.  {\em Scalar curvature and stability of toric
varieties}, J. Differential Geom. {\bf 62} (2002), no. 2, 289--349
\bibitem[D2] {D2} Donaldson, S.K.
{\em Lower bounds for the Calabi energy},  J. Differential Geom.  {\bf 70}  (2005),  no. 3, 453--472
\bibitem[F]{F} Futaki, A. {\em K\"ahler-Einstein metrics and integral invariants}, Lecture Notes in
Mathematics, 1314. Springer-Verlag, Berlin, 1988
\bibitem[H] {H} Hamilton, R.S., {\em The Ricci flow on surfaces},
Contemp. Math. {\bf 71} (1988) 237--261
\bibitem[L]{L} Li, H. {\em On the lower bound of the K energy and F functional}, preprint,
arXiv:math.DG/0609725
\bibitem[Pa] {Pa} Pali, N. 
{\em Characterization of Einstein-Fano manifolds via
the K\"ahler-Ricci flow}, arXiv: math.DG/0607581.
\bibitem[P1]{P1} Perelman, G. {\em The entropy formula for the Ricci flow and its geometric applications},
preprint, arXiv:math.DG/0211159
\bibitem[P2]{P2} Perelman, G. unpublished work on the K\"ahler-Ricci flow
\bibitem[PSSW]{PSSW} Phong. D.H., Song. J., Sturm, J. and Weinkove, B.
{\em The Moser-Trudinger inequality on K\"ahler-Einstein manifolds},
arXiv: math.DG/0604076, to appear in Amer. J. Math.
\bibitem[PSS] {PSS} Phong, D.H., Sesum, N. and Sturm, J.,
{\em Multiplier ideal sheaves and the K\"ahler-Ricci flow},
arXiv: math.DG/0611794, to appear in Comm. Anal. and Geometry
\bibitem[PS]{PS} Phong, D.H. and Sturm, J.  {\em On stability and the convergence of the K\"ahler-Ricci flow},  J. Differential Geom.  {\bf 72} (2006),  no. 1, 149--168
\bibitem[ST]{ST} Sesum, N. and Tian, G. {\em Bounding scalar curvature and diameter along the K\"ahler-Ricci
flow (after Perelman) and some applications}, preprint
\bibitem[T]{T} Tian, G. {\em K\"ahler-Einstein metrics with positive scalar curvature},
Invent. Math. {\bf 130} (1997), 1--37
\bibitem[TZ1]{TZ1}
Tian, G. and Zhu, X. {\em A nonlinear inequality of Moser-Trudinger
type}, Calc. Var. {\bf 10} (2000), 349--354
\bibitem[TZ2]{TZ2} Tian, G. and Zhu, X. {\em Convergence of K\"ahler-Ricci flow}, J. Amer. Math. Soc. {\bf 20} (2007), no. 3, 675--699
\bibitem[To]{To} Topping, P. {\em Lectures on the Ricci flow}, London Mathematical Society Lecture Note
Series, 325. Cambridge University Press, Cambridge, 2006
\bibitem[WZ] {WZ} Wang, X.J. and Zhu, X.
{\em K\"ahler-Ricci solitons on toric manifolds with positive first Chern class},
Advances Math. {\bf 188} (2004) 87--103
\bibitem[Y1]{Y1} Yau, S.-T. {\em On the Ricci curvature of a compact K\"ahler
manifold and the complex Monge-Amp\`ere equation, I}, Comm. Pure
Appl. Math. {\bf 31} (1978), 339--411
\bibitem[Y2]{Y2} Yau, S.-T. {\em Open problems in geometry}, Proc. Symposia Pure
Math. {\bf 54} (1993), 1--28 (problem 65)
\bibitem[Z]{Z} Zhu, X. {\em K\"ahler-Ricci flow on a toric manifold with positive first Chern class},
preprint,  arXiv:math.DG/0703486

\v\v
$^*$ Department of Mathematics\\
Columbia University, New York, NY 10027\\


$^{**}$ Department of Mathematics \\
Johns Hopkins University, Baltimore, MD 21218\\



$^{\dagger}$ Department of Mathematics \\
Rutgers University, Newark, NJ 07102\\


$^{\ddagger}$ Department of Mathematics \\
Harvard University, Cambridge, MA 02138




\end{thebibliography}
\end{document}